\documentclass[12pt]{article} 

\usepackage{latexsym}
\usepackage{amssymb}
\usepackage{euscript}
\usepackage[dvips]{graphics}
\usepackage{epsf}

\let\cal=\mathcal      

\def\mcc{M\raise.5ex\hbox{c}C}
\def\mccarthy{M\raise.5ex\hbox{c}Carthy}
\def\Hu{\H}
\def\sz{Szeg\Hu{o} }

\def\eg{{\it e.g. }}
\def\ie{{\it i.e. }}

\def\atn{{\A}^2(\nu)}
\def\htm{H^2(\mu)}

\def\ltm{L^2(\mu)}
\def\h{{\cal H}}
\def\hk{{\cal H}_k}

\def\N{{\cal N}}

\def\Om{\Omega}

\def\l{\lambda}
\def\z{\zeta}

\def\vare{\varepsilon}

\let\ii=\i
\let\i=\infty

\def\la{\langle}
\def\ra{\rangle}
\def\={ = }

\def\A{{\cal A}}

\def\C{\mathbb C}

\def\T{\mathbb T}
\def\D{\mathbb D}

\def\be{\setcounter{equation}{\value{theorem}} \begin{equation}}
\def\ee{\end{equation} \addtocounter{theorem}{1}}
\def\beq{\begin{eqnarray*}}
\def\eeq{\end{eqnarray*}}
\def\se{\setcounter{equation}{\value{theorem}}} 
\def\att{\addtocounter{theorem}{1}}

\def\bs{}

\def\quest{\bs \att {{\bf Question \thetheorem \ }} }

\def\bp{{\sc Proof: }}
\def\ep{{}{\hfill $\Box$} \vskip 5pt \par}

\def\obc{{\sc Proof of Claim: }}
\def\oec{{}{\hspace*{\fill} $\lhd$} \vskip 5pt \par}
\def\bl{\begin{lemma}}
\def\el{\end{lemma}}
\def\bt{\begin{theorem}}
\def\et{\end{theorem}}
\def\bprop{\begin{prop}}
\def\eprop{\end{prop}}
\def\bd{\begin{definition}}
\def\ed{\end{definition}}
\def\br{\begin{remark}}
\def\er{\end{remark}}
\def\bexer{\begin{exercise}}
\def\eexer{\end{exercise}}

\newtheorem{theorem}{Theorem}[section]
\newtheorem{prop}[theorem]{Proposition}
\newtheorem{lemma}[theorem]{Lemma}

\newtheorem{definition}[theorem]{Definition}

\newcommand{\U}{\cal U}
\newcommand{\g}{\gamma}
\newcommand{\lN}{\lambda_{N+1}}
\newcommand{\wN}{w_{N+1}}
\newcommand{\hib}{H^\infty (\D^2)}

\begin{document}
\setlength{\baselineskip}{21pt}
\title{ Distinguished Varieties}
\author{Jim Agler
\thanks{Partially supported by the National Science Foundation}\\
U.C. San Diego\\
La Jolla, California 92093
\and
John E. M\raise.5ex\hbox{c}Carthy
\thanks{Partially supported by National Science Foundation Grant
DMS 0070639}\\
Washington University\\
St. Louis, Missouri 63130}
\date{August 20 2004}

\bibliographystyle{plain}

\maketitle
\begin{abstract}
A distinguished variety is a variety that exits the bidisk through
the distinguished boundary. We show that And\^o's inequality for
commuting matrix contractions can be sharpened to looking at the
maximum modulus on a distinguished variety, not the whole bidisk.
We show that uniqueness sets for extremal Pick problems on the
bidisk always contain a distinguished variety.
\end{abstract}

\baselineskip = 18pt

\setcounter{section}{-1}
\section{Introduction}

In this paper, we shall be looking at a special class of bordered
(algebraic) varieties that are contained in the bidisk $\D^2$ in $\C^2$.
\bd A non-empty 
set $V$ in $\C^2$ is a {\em distinguished variety} 
if there is a polynomial $p$ in $\C[z,w]$ 
such that 
$$
V \= \{ (z,w)  \in  \D^2 \ : \ p(z,w) = 0 \} 
$$
and such that
\be
 \overline{V} \cap \partial ( \D^2) \=  \overline{V} \cap (\partial \D)^2
.
\label{eqa1}
\ee
\ed
Condition (\ref{eqa1}) means that the variety exits the bidisk
through the distinguished boundary of the bidisk, the torus.
We shall use $\partial V$ to denote the set given by (\ref{eqa1}):
topologically, it is the boundary of $V$ within $Z_p$, the zero set of
$p$, rather than in all of $\C^2$. We shall always assume that $p$
is chosen to be minimal, \ie so that no irreducible component of
$Z_p$ is disjoint from $\D^2$ and so that $p$ has no repeated
irreducible factors.
\bs
Why should one single out distinguished varieties
from other bordered varieties?

One of the most important results in operator theory is 
T.~And\^o's inequality \cite{and63} (see also \cite{colwer99} and
\cite{szn-foi}). This says that if $T_1$ and $T_2$ are commuting
operators, and both of them are of norm $1$ or less, then for any
polynomial $p$ in two variables, the inequality
\be
\label{eqa2}
\| p(T_1,T_2) \| \ \leq \ \| p \|_{\D^2}  
\ee
holds. And\^o's inequality is essentially
equivalent to the commutant lifting theorem of B.~Sz.-Nagy and
C.~Foia\c{s} \cite{sznfoi68a} --- see \eg \cite{niva98} for a
discussion of this. 

Our first main result, Theorem~\ref{thmc1}, 
is that if $T_1$ and $T_2$ are matrices, then the
inequality~(\ref{eqa2})
can be improved to 
$$
\| p(T_1,T_2) \| \ \leq \ \| p \|_V ,
$$
where $V$ is some distinguished variety depending on
$T_1$ and $T_2$.
Indeed, in the proof of the theorem, we construct co-isometric
extensions of the matrices that naturally live on this distinguished
variety. So when studying bivariable matrix theory, rather than 
operator theory, one is led inexorably to study distinguished
varieties.

Conversely, in Theorem~\ref{thmz1}, we show that all distinguished varieties
can be represented as
$$ \{ (z,w) \in \D^2 \ : \
\det(\Psi(z) - w I) = 0 \}
$$
for some analytic matrix-valued function $\Psi$
on the disk that is unitary on $\partial \D$. 
This shows that the study of distinguished varieties leads back to
operator theory.
\bs
Consider the natural notion of isomorphism of two distinguished
varieties, namely that there is a
biholomorphic bijection between them.
\bd
A function $\Phi$ is holomorphic on a set $V$ in $\C^2$
if, at every point $\l$ in $V$, there is a non-empty
ball $B(\l,\vare)$ centered at $\l$
and an analytic mapping of two variables
defined on $B(\l,\vare)$ that agrees with
$\Phi$ on $B(\l,\vare) \cap V$.
\ed 
\bd Two
distinguished varieties $V_1$ and $V_2$ are {\em isomorphic} if there is
a function $\Phi$ that is holomorphic on $V_1$ and continuous
on $\overline{V_1}$ such that $\Phi$ is a bijection from 
$\overline{V_1}$ onto
$\overline{V_2}$ 
and such that $\Phi^{-1}$ is holomorphic on $V_2$.
\label{defa2}
\ed
(The requirement that $\Phi^{-1}$ be holomorphic does not 
follow automatically from the holomorphicity of $\Phi$ ---
consider \eg $V_1 \= \{ (z,z)  :  z  \in  \D \}$
and $V_2 \= \{ (z^2,z^3)  :  z  \in  \D \}$, which are not
isomorphic.)

By the maximum modulus principle, $\Phi$ must map the boundary of
$V_1$ onto the boundary of $V_2$. It follows that $\Phi =
(\phi_1,\phi_2)$ is a pair of {\em inner functions}, \ie a pair of
holomorphic scalar-valued functions that each have modulus 
one on $\partial V_1$. So studying isomorphism classes of
distinguished varieties is closely connected to the rich structure
of inner functions. 

W.~Rudin has studied when an arbitrary finite
Riemann surface $R$ is isomorphic to a distinguished variety, in
the sense that there is an unramified pair of separating inner
functions on $R$ that are continuous on $\overline{R}$
\cite{rud69b}. His results show, for example, that a finitely
connected planar domain is isomorphic to a distinguished variety
if and only if the domain is either a disk or an annulus.
He also showed that for every $n \geq 1$, there is a finite Riemann
surface $R$ that is topologically an $n$-holed torus minus one disk,
and such that $R$ is isomorphic to a distinguished variety.

In Section~\ref{secb} we show that, under fairly general
conditions, a pair of ``inner'' functions $(\phi_1,\phi_2)$ 
on a set $X$ must map $X$
into a distinguished variety (\ie the algebraic relation on the
$\phi_i$'s comes for free). 

\bs
A third reason to study distinguished varieties comes from
considering the Pick problem on the bidisk. 
This is the problem of deciding, given points $\l_1, \dots , \l_N$
in $\D^2$, and values $w_1, \dots, w_N$ in $\C$, whether there is a
function in $\hib$, the bounded analytic functions on $\D^2$,
that interpolates the data and is of norm at most one. The problem
is called {\it extremal} if there is an interpolating function of norm
exactly one, but not less.

If an extremal Pick
problem is given, the solution may or may not be unique (see
Section~\ref{secd} for an example).
Our second main result is
Theorem~\ref{thmd1}, where we show that 
there is always a distinguished variety on which the solution is
unique. 

One can think then of the Pick problem
as having two parts: 

\noindent
(a) Solve the problem on the distinguished
variety where the solution is unique.

\noindent
(b) Parametrize all the extensions 
of the solution to the whole bidisk.

We give a formula (\ref{eqd8}) for Problem (a).
The extension problem (b) is non-trivial: unless the distinguished
variety is isomorphic to a disk, there will always be some functions that
cannot be extended to the whole bidisk without increasing the norm
\cite{agmc_vn}. Obviously a function arising from a Pick problem
will be extendable, but what distinguishes such functions remains
mysterious.

If one starts with an inner function on $V$ and wants to extend this
to a rational inner function on $\D^2$, there can be more than one
extension. However, there is a restriction
on the degree, given by Theorem~\ref{thmb2}. If the variety is of rank
$(n_1,n_2)$, \ie there are generically $n_1$ sheets above every first
coordinate and $n_2$ above every second coordinate, then 
any regular rational inner extension of degree $(d_1,d_2)$
must have
$d_1 n_1 + d_2 n_2$
equal to the number of zeroes the original function had on $V$.

An {\it admissible kernel $K$} on a set $\{\l_1, \dots, \l_N\}$ in
$\D^2$ is
an $N$-by-$N$ positive definite matrix 
such that
$$
[ ( 1 - \l_i^r \overline{\l_j^r}) K_{ij} ] \ \geq \ 0 \qquad r =
1,2 .
$$
It is known \cite{ag1,ampi} that studying all the admissible
kernels on a set is essential to understanding the Pick problem.
A key idea in the proof of Theorem~\ref{thmd1} is that 
every admissible kernel automatically extends to a distinguished
variety.

Distinguished varieties have been studied in a somewhat more
abstract and general setting by J.~Ball and V.~Vinnikov
\cite{bv00}.
They have a determinental representation that is analogous to
Theorem~\ref{thmz1}.

We would like to thank the referees for many
valuable remarks.

\section{Representing Distinguished Varieties}
\label{secz}

Let $V$ be a distinguished variety. We say a function $f$ is {\em
holomorphic on $V$} if, for every point of $V$, there is an
open ball $B$ in $\C^2$ containing the point, and a holomorphic
function $\phi$ of two variables on $B$, such that $\phi |_{B \cap
V} = f |_{B \cap V}$. We shall use $A(V)$ to denote the Banach
algebra of functions that are holomorphic on $V$ and continuous on
$\overline{V}$. 
This is a uniform algebra on $\partial V$, \ie a closed unital
subalgebra of $C(\partial V)$ that separates points.
The maximal ideal space of $A(V)$ is $\overline{V}$. 

If $\mu$ is a finite measure on a distinguished variety $V$, let $\htm$
denote the closure in $\ltm$ of the polynomials.
If $\Omega$ is an open subset of a Riemann surface $S$, and $\nu$ is a
finite
measure on $\overline{\Omega}$, let $\atn$ denote the closure in
$L^2(\nu)$ of $A(\Om)$, 
the functions that are  holomorphic on $\Omega$ and
continuous on $\overline{\Omega}$. We say a point $\l$ is a bounded
point evaluation for $\htm$ (or $\atn$) if evaluation at $\l$,
{\em a priori} defined only for a dense set of
analytic functions, extends continuously to the whole Hilbert
space. If $\l$ is a bounded point evaluation, we call the function
$k_\l$ that has the property that 
$$
\la f, k_\l \ra \= f(\l) 
$$
the {\em evaluation functional at $\l$}. 



The following lemma is well-known. 
It is valid in much greater generality, but this will
suffice for our purposes.
If the boundary of $\Omega$
consists of closed analytic curves, the lemma follows from
J.~Wermer's proof \cite{wer64} that $A(\Omega)$ is hypo-Dirichlet,
and the description of representing measures for hypo-Dirichlet
algebras given by P.~Ahern and D.~Sarason in \cite{ahesar1}.
(Actually Wermer's proof extends without difficulty to the case
where the boundary is just piecewise $C^2$, but we shall not need
this fact).
For a detailed description of the measures in this case, see
K.~Clancey's paper \cite{cl91}. 

\bl
\label{lemz00}
Let $S$ be a compact Riemann surface. Let $\Om \subseteq S$
be a domain 
whose boundary is a finite union of piecewise smooth Jordan curves.
Then there exists a measure $\nu$ on $\partial \Om$ 
such that
every $\l$ in $\Om$ is a bounded point evaluation for $\atn$, and such that the linear
span of the evaluation functionals is dense in $\atn$.
\el
\bp
%
%
Because its boundary is nice, $\Om$ is regular for the
Dirichlet problem (see \eg \cite[Section IV.2]{farkra}).
Let $\nu$ be harmonic measure for $\Om$ with respect 
to some fixed base-point.
Then by Harnack's inequality, harmonic measure for any other 
point in the domain is
boundedly absolutely continuous with respect to $\nu$. As harmonic evaluation
functionals are {\em a fortiori} analytic evaluation functionals, we get that
every point of $\Om$ is a bounded point evaluation (with an $L^\i$ evaluation
functional) for $\atn$.

Ahern and Sarason \cite[p.159]{ahesar1} proved that the 
span of the evaluation functionals is dense.
Their argument, in brief, was to 
find an exhaustion $\Om_n$ of $\Om$, \ie an increasing family of
open sets, each contained compactly in the next, whose union was
$\Om$. Let $\nu_n$ be harmonic measure for each $\Om_n$, with
respect to the same fixed base-point. Then they showed that 
for every $u$ in $L^1(\partial \Om, \nu)$, its norm was equal to
$$
\lim_{n \to \i} \int |\hat{u}| d \nu_n ,
$$
where $\hat{u}$ is the harmonic extension of $u$ to $\Om$.
In particular, any function in $\atn$ that vanishes identically on
$\Om$ must be the zero function.

\ep

\bl
\label{lemz0}
Let $V$ be a distinguished variety.
There is a measure $\mu$ on $\partial V$ 
such that every point in $V$ 
is a bounded point evaluation for $\htm$, and such that the
span of the evaluation functionals is dense in $\htm$.
\el
\bp
Let $p$ be the minimal polynomial such that
$V$ is the intersection of $Z_p$
with $\D^2$. Let $C$ be the projective closure of $Z_p$
in $\C {\mathbb P}^2$. Let $S$ be the desingularization of $C$.
This means $S$ is a compact Riemann surface (not connected if $C$
is not irreducible) and there is a holomorphic function
$\phi : S \to C$ that
is biholomorphic from $S'$ onto $C'$ and
finite-to-one from $S \setminus S'$ onto $ C \setminus C'$.
Here $C'$ is the set of non-singular points in $C$, and $S'$ is the
preimage of $C'$. See \eg \cite{fis01} or \cite{gh78} for details
of the desingularization.


Let $\Om = \phi^{-1} ( V)$. Then $\partial \Om$ is a finite union
of disjoint curves, each of which is analytic except possibly at a
finite number of cusps. 
Let $\nu$ be the measure from
Lemma~\ref{lemz00} (or the sum of these if $\Om$ is not connected).

The desired measure $\mu$ is the push-forward of $\nu$ by $\phi$,
\ie it is defined by $\mu(E) = \nu (\phi^{-1}(E))$.
Indeed, if $\l$ is in $V$ and $\phi(\z) = \l$, let
$k_\z \nu$ be a representing measure for $\z$ in $A(\Om)$.
Then the function $k_\z \circ \phi^{-1}$ is defined 
$\mu$-a.e., and satisfies
$$
\int_{\partial V} p (k_\z \circ \phi^{-1})  d \mu
\=\
\int_{\partial \Om} (p \circ \phi)  k_\z d\nu
\= p\circ\phi(\z) \= p(\l) .
$$ 
\ep
Note that $\{ g \circ \phi \ : \ g  \in  A(V) \}$ is a finite
codimensional subalgebra of $A(\Om)$. For a description of what finite
codimensional subalgebras look like, see Gamelin's paper \cite{gam68}.

For positive integers $m$ and $n$, let
\be
\label{eqz05}
U \= \left( \begin{array}{cc} A & B \\
C & D \end{array} \right) \ : \ \C^m \oplus \C^n
\ \to \ \C^m \oplus \C^n
\ee
be an $(m+n)$-by-$(m+n)$ unitary matrix.
Let
\be
\label{eqz1}
\Psi(z) =  A + zB(I- zD)^{-1}C
\ee
be the $m$-by-$m$ matrix valued
function defined on the unit disk $\D$ by the entries of $U$.
This is called the {\it transfer function} of $U$.
Because $U^\ast U = I$, a calculation yields
\be
\label{eqz2}
I - \Psi(z)^\ast \Psi (z) \=
(1 - |z|^2)\
C^\ast ( I - \bar z D^\ast)^{-1}
(I-z D)^{-1}C ,
\ee
so $\Psi(z)$ is a rational matrix-valued function that is unitary
on the unit circle and contractive on the unit disk. 
Such functions are called rational matrix
inner functions, and it is well-known that all 
rational matrix inner functions 
have the form (\ref{eqz1}) for some unitary matrix
decomposed as in (\ref{eqz05}) --- see \eg \cite{ampi} for a proof.

Let $V$ be the set
\be
\label{eqz4}
V \= \{ (z,w) \in \D^2 \ : \
\det(\Psi(z) - w I) = 0 \}.
\ee
We shall show that $V$ is a distinguished variety, and that every distinguished
variety arises this way.

\bl
\label{lemz1}
Let
$$
U' \= \left( \begin{array}{cc} D^\ast & B^\ast \\
C^\ast & A^\ast \end{array} \right) \ : \ \C^n \oplus \C^m
\ \to \ \C^n \oplus \C^m,
$$
let
$$
\Psi'(z) \=  D^\ast + zB^\ast(I- z A^\ast)^{-1}C^\ast,
$$
and let
$$
V' \=
\{ (z,w) \in \D^2 \ : \
\det(\Psi'(w) - z I) = 0 \}.
$$
Then $V = V'$.
\el
\bp
The point $(z,w)  \in  \D^2$
is in $V$ iff there is a non-zero vector $v_1$ in
$\C^m$ such that
\be
\label{eqb5}
\left[ A + zB (1 - zD)^{-1} C \right] v_1 \= w v_1 .
\ee
Claim: (\ref{eqb5}) holds if and only if there is a non-zero vector
$v_2$ in $\C^n$ such that
\be
\label{eqb6}
\left( \begin{array}{cc} A & B \\
C & D \end{array} \right)
\left( \begin{array}{c} v_1 \\
z\ v_2 \end{array} \right)
\=
\left( \begin{array}{c} w\ v_1 \\
v_2 \end{array} \right) .
\ee

\obc If (\ref{eqb6}) holds, then solving gives
(\ref{eqb5}). Conversely, if (\ref{eqb5}) holds, define $$
v_2 \= (I - zD)^{-1} C v_1.
$$
Then (\ref{eqb6}) holds. Moreover, if $v_2$ were
$0$, then $v_1$ would be in the kernel of $C$ and be a
$w$-eigenvector of $A$.
As $A^\ast A + C^\ast C = I$, this would force $|w| =1$,
contradicting the fact that
$(z,w)  \in  \D^2$.
\oec

Given the claim, the point $(z,w)$ is in $V'$ iff there are
non-zero vectors $v_1 $ and $v_2$ such that
\be
\label{eqb65}
\left( \begin{array}{cc} D^\ast & B^\ast \\
C^\ast & A^\ast \end{array} \right)
\left( \begin{array}{c} v_2 \\
w\ v_1 \end{array} \right)
\=
\left( \begin{array}{c} z\ v_2 \\
v_1 \end{array} \right) .
\ee
Interchanging coordinates, (\ref{eqb65}) becomes
\be
\label{eqb7}
\left( \begin{array}{cc} A^\ast & C^\ast \\
B^\ast & D^\ast \end{array} \right)
\left( \begin{array}{c} w\ v_1 \\
v_2 \end{array} \right)
\=
\left( \begin{array}{c} v_1 \\
z \ v_2 \end{array} \right) .
\ee
Clearly, (\ref{eqb6}) and (\ref{eqb7})
are equivalent.
\ep

Note that if $C$ has a non-trivial kernel $\N$, then 
(\ref{eqz2}) shows that $\Psi(z)$ is isometric on $\N$ for all $z$,
so by the maximum principle is equal to a constant isometry with
initial space $\N$. If $C$ has a trivial kernel, we say $\Psi$ is
{\em pure}. Every rational inner function decomposes into the
direct sum of a pure rational inner function and a unitary matrix
--- see \eg \cite{szn-foi}.
Since $A^\ast A + C^\ast C = I$, we see that $C$ has no kernel iff
$\| A \| < 1$. Since $A A^\ast + B B^\ast = I$, this in turn is
equivalent to $B^\ast$ having no kernel.
Therefore $\Psi$ is pure iff $\Psi'$ is.

\bt
\label{thmz1}
The set $V$, defined by (\ref{eqz4}) for some rational matrix inner
function $\Psi$, is a distinguished variety. Moreover, every
distinguished variety can be represented in this form.
\et
\bp
Suppose  $V$ is given by (\ref{eqz4}), and 
that $(z,w)$ is in $\overline{V}$.
Without loss of generality, we can assume that $\Psi$ is pure.
Indeed, any unitary summand of $\Psi$ would add sheets to the 
variety $\det (\Psi(z) - wI) = 0$ of the type $\C \times \{ w_0
\}$, for some unimodular $w_0$. These sheets are all disjoint from
the open bidisk $\D^2$.

If $|z| < 1$, equation~(\ref{eqz2}) then shows that $\Psi(z)$ is a
strict contraction, so all its eigenvalues must have modulus
less than $1$, and so  $|w| < 1 $ also. To prove that $|w| < 1$ implies
$|z| < 1$, just apply the same argument to $V'$.
Therefore (\ref{eqa1}) holds, and $V$ is a distinguished variety.

To prove that all distinguished varieties arise in this way, let
$V$ be a distinguished variety. Let $\mu$ be the measure from
Lemma~\ref{lemz0},
and let $\htm$ be the closure of the
polynomials in $L^2(\mu)$.
The set of bounded point evaluations for $\htm$ is precisely $V$.
(It cannot be larger, because $\overline{V}$ is polynomially convex,
and Lemma~\ref{lemz0} ensures that it is not smaller).

Let $T = (T_1,T_2)$ be the pair of operators on $\htm$ given by
multiplication by the coordinate functions. They are pure commuting
isometries\footnote{
A pure isometry $S$ is one that has no unitary summand; this is the
same as requiring that $\cap_{i=1}^\i {\rm ran}(S^i )  =  \{ 0
\}.$
} because the span of the evaluation functionals is
dense. The joint eigenfunctions of their adjoints are the
evaluation functionals.

By the Sz.-Nagy-Foia\c{s} model theory \cite{szn-foi}, $T_1$ can
be modelled as $M_z$, multiplication by the independent variable
$z$ on $H^2 \otimes \C^m$, a vector-valued Hardy
space on the unit circle. In this model, $T_2$ can be modelled as $M_\Psi$,
multiplication by $\Psi(z)$ for some pure rational matrix 
inner function $\Psi$. 
A point $(z,w)$ in $\D^2$ is a bounded point evaluation for $\htm$
iff $(\bar z, \bar w)$ is a joint eigenvalue for
$(T^\ast_1,T^\ast_2)$. In terms of the unitarily equivalent
Sz.-Nagy-Foia\c{s} model, this is equivalent to $\bar w$ being an
eigenvalue of $\Psi(z)^\ast$.

Therefore $$
V = \{ (z,w)  \in  \D^2 \ : \ \det (\Psi(z) - wI)  =  0 \},
$$
as desired.
\ep
If $\Psi$ is the transfer function of a unitary $U$ as in
(\ref{eqz05}), and $\Psi$ is pure, 
we shall say that $V$ is of {\em rank $(m,n)$}. This
means that generically there are $m$ sheets above each $z$, and $n$
sheets above each $w$.

\section{
Inner Functions}
\label{secb}
Rudin's results \cite{rud69b} show that planar annuli can be mapped
isomorphically into distinguished varieties by a pair of inner
functions. The advantage of doing
this is that the coordinate functions are then easier to deal with
than the original inner functions. 
Inner functions on a finite bordered Riemann surface can be shown
to satisfy an algebraic equation. In this section, we show that
even without the Riemann surface structure, inner functions must
satisfy an algebraic equation. The result is reminiscent
of Livsic's Cayley-Hamilton theorem for a pair of commuting
operators with finite rank imaginary parts --- see \eg the book
\cite{lkmv}.


Let $X$ be a set.
By a {\em kernel} on $X$ we mean a self-adjoint map
$k: X \times X \to \C$ that is positive definite, in the sense that
for any finite set $\{\l_1,\dots, \l_N\}$ of distinct points
in $X$, the self-adjoint matrix
$k(\l_j,\l_i)$ is positive definite. 
Given any kernel $k$, there is a Hilbert space $\h_k$ of functions
on $X$ for which $k$ is the reproducing kernel, \ie
$$
\la f(\cdot), k(\cdot, \l) \ra
\= f(\l) \qquad \forall  f  \in \h_k,\ \forall  \l  \in X .
$$
(For details of the passage between a kernel and a Hilbert function
space, see \eg \cite{ampi}).

Let $\phi_1$ and $\phi_2$ be functions on $X$ with modulus less
than one at every point. Assume that we can
find some kernel $k$ on $X$ so that multiplication by each $\phi_i$
is a pure
isometry
on $\hk$ with finite dimensional cokernel. For example, $X$
could be a distinguished variety, the $\phi_i$'s could be the
coordinate functions, and $\hk$ could be the closure of the
polynomials in $L^2(\partial X)$.  Or, $X$ could be a 
smoothly bounded planar domain, the $\phi_i$'s could 
be inner functions that are
continuous on $\overline{X}$ and have finitely many zeroes, and
$\hk$ could be the closure in $L^2(\partial X)$ of the rational
functions with poles off $\overline{X}$.

Let $e_1, \dots, e_m$ be an orthonormal basis for $(\phi_1 \hk
)^\perp$.  Then 
$$
\{ \phi_1^i e_j \ : \ i  \in  {\mathbb N},\ 1  \leq  j 
\leq  m \}
$$
is an orthonormal basis for $\hk$.
So by Bergman's formula \cite[Prop 2.18]{ampi}, 
\se\att
\begin{eqnarray}
\nonumber
k(\z,\l) &\=& \sum_{i=0}^\i \sum_{j=1}^m \phi_1^i(\z)  e_j (\z)
\overline{\phi_1^i(\l)  e_j (\l)} \\
&=& \frac{\sum_{j=1}^m  e_j (\z) \overline{  e_j (\l)}}
{1-\phi_1(\z) \overline{\phi_1(\l)}}.
\label{eqb1}
\end{eqnarray}

Similarly, if $f_1, \dots, f_n$ is an orthonormal basis for
$(\phi_2 \hk)^\perp$, we get
\be
\label{eqb2}
k(\z,\l) \= \frac{\sum_{j=1}^n  f_j (\z) \overline{  f_j (\l)}}
{1-\phi_2(\z) \overline{\phi_2(\l)}}.
\ee
Equating the right-hand sides of (\ref{eqb1}) and (\ref{eqb2}) and
cross-multiplying, we get
\se\att
\begin{eqnarray}
\lefteqn{
\sum_{j=1}^m e_j(\z) \overline{e_j(\l)} 
\ + \ 
\sum_{i=1}^n \phi_1(\z) f_i(\z) \overline{ \phi_1 (\l) f_i(\l)}
\nonumber} \\
&&\=
\sum_{j=1}^m \phi_2(\z) e_j(\z) \overline{\phi_2(\l) e_j(\l)} 
\ + \ 
\sum_{i=1}^n  f_i(\z) \overline{  f_i(\l)} .
\label{eqb3}
\end{eqnarray}

Let $f(\z)$ be the vector in $\C^n$ with components $f_1(\zeta),
\dots , f_n(\z)$, and let $e(\z) = (e_1(\z), \dots, e_m(\z))^t$.
Then (\ref{eqb3}) can be rewritten as saying that the map
\beq
U \= \left( \begin{array}{cc} A & B \\
C & D \end{array} \right) \ : \ \C^m \oplus \C^n 
&\ \to \ &\C^m \oplus \C^n \\
\left( \begin{array}{c} 
e(\z) \\
\phi_1 (\z) f(\z) 
\end{array} \right)
&\mapsto &
\left( \begin{array}{c} \phi_2 (\z) e(\z) \\
f(\z) \end{array} \right)
\eeq
is an isometry on the linear span of the vectors
$$
\left\{ \left( \begin{array}{c} e(\z) \\ \phi_1 (\z) f(\z) 
 \end{array} \right) \ : \ \z \in X \right\} .
$$
Even if these vectors do not span all of $\C^m \oplus \C^n$, we can
always extend $U$ to be a unitary from $\C^m \oplus \C^n$ onto $\C^m
\oplus \C^n$, and we shall assume that we have done this.


Let 
\be
\label{eqb37}
\Psi(z) =  A + zB(I- zD)^{-1}C
\ee
be the $m$-by-$m$ matrix valued
function defined on the unit disk $\D$ that is the transfer
function of $U$.
Moreover, we have
$$
\Psi(\phi_1(\z) )  e(\z) \= \phi_2(\z)  e(\z) .
$$
Therefore the points $(\phi_1(\z), \phi_2(\z) )$
all lie in the set
\be
\label{eqb45}
V \= \{ (z,w) \in \D^2 \ : \
\det(\Psi(z) - w I) = 0 \},
\ee
which we know from Theorem~\ref{thmz1} is a distinguished variety.
Thus we have proved:
\bt
\label{thmb1}
Let $\hk$ be a reproducing kernel Hilbert space on a set $X$. Let 
$\phi_1$ and $\phi_2$ be multipliers of $\hk$ such that 
multiplication by each $\phi_i$ is a pure isometry with finite
dimensional cokernel, and such that $|\phi_i (\z)| < 1\quad \forall
 \z  \in  X$.
With notation as above, the function
$$
\zeta \ \mapsto \ (\phi_1(\z), \phi_2(\z))
$$
maps $X$ into the distinguished variety $V$ given by
(\ref{eqb45}).
\et

\bs
Note that applying Theorem~\ref{thmb1} to $\htm$, the space in 
Lemma~\ref{lemz0}, we get the second part of Theorem~\ref{thmz1}.

If $V$ is a distinguished variety, an inner function on $V$
may or may not extend to an inner function on $\D^2$. If it does
extend, the extension may not be unique. It is curious, however,
that there is a rigidity in the degree of this extension.
Let $\phi$ be a rational inner function on $\D^2$. Then it can be
represented as
\be
\label{eqb36}
\phi(\z) \= \frac{\z^d \overline{p(\frac{1}{\bar \z})}}{p(\z)}
\ee
for some polynomial $p$ that does not vanish on $\D^2$
\cite{rud69}, where $\z = (\z_1,\z_2)$ and $d$ is a multi-index.
The representaion is not unique --- \eg taking $p(z,w) = i (z^2 -
w^2)$ and $d = (2,2)$, one gets the constant function $1$.
The representation will be unique if $p$ is restricted so that
$Z_p \cap \T^2$ is finite. In this event,
we shall call $d = (d_1,d_2)$ the degree of $\phi$.

If $\phi$ is an inner function in $A({\D^2})$,
then it is rational and moreover the function
$p$ will not vanish on $\overline{\D^2}$ \cite[Thm. 5.2.5]{rud69};
we shall call such a  function {\it regular}.
\bt
\label{thmb2}
Let $V$ be a variety of rank $n = (n_1,n_2)$, and let $\phi$ be a
regular
rational
inner function on $\D^2$ of degree $d$. Then $\phi$ restricted to
$V$ has exactly $n \cdot d =  n_1 d_1 + n_2 d_2$
zeroes, counting multiplicities.
\et
\bp
By applying an automorphism of $\D^2$, we can assume that $(0,0)$ is not in
$V$ and that all points with first or second coordinate $0$ are regular.

Consider first the case $\phi(z,w) = z^{d_1} w^{d_2}$, \ie 
$p \equiv 1$ in (\ref{eqb36}). Then at each of the $n_1$ points in
$V$ with second coordinate $0$ has a zero of multiplicity $d_1$, and each of the 
$n_2$ points in $V$ with first coordinate $0$ has a zero of multiplicity $d_2$.

Now let $p$ be an arbitrary polynomial that does not vanish on $\overline{\D^2}$,
normalized so that $p(0,0)=1$.
Let $p_r(\z) = p(r \z)$, and
$$
\phi_r(\z) \= \frac{\z^d \overline{p_r(\frac{1}{\bar \z})}}{p_r(\z)}
.$$
As $r$ increases from $0$ to $1$, the function $\phi_r$ changes continuously from
$\z^d$ to $\phi$. As each $\phi_r$ is in $A(V)$ and is inner, the number of zeroes 
must remain constant.
\ep

Example.
Let $V$ be the distinguished variety $\{ z^2 = w^3 \}$, of rank 
$(3,2)$. The inner function $\phi(z,w) = z^2$ can be extended to either the function
$z^2$ of 
degree $(2,0)$ or $w^3$ of degree $(0,3)$. In either event, $n \cdot d = 6$.

\section{A sharpening of And\^o's inequality}
\label{secc}

\bt
\label{thmc1}
Let $T_1 $ and $T_2$ be commuting contractive matrices, neither of
which has eigenvalues of modulus $1$.
Then there is a distinguished variety $V$ such that, for any
polynomial $p$ in two variables, the inequality
\be
\label{eqc1}
\| p(T_1, T_2) \| \ \leq \ \| p \|_V 
\ee
holds.
\et
\bp
Let the dimension of the space on which the matrices act be $N$.

(i) First, let us assume that each $T_r$ has $N$
linearly independent unit eigenvectors, $\{ v_j \}_{j=1}^N$.
So we have
$$
T_r v_j \= \lambda_j^r v_j, \qquad  r=1,2 \quad 1\leq j \leq N ,
$$
for some set of scalars $\{ \l_j^r \}$.
As each $T_r$ is a contraction, we have $I -  T_r^\ast T_r $ is
positive semidefinite, so
\be
\label{eqc2}
\la ( I -  T_r^\ast T_r ) v_j, v_i \ra 
\= (1 -\overline{\l_i^r} \l_j^r ) \la v_j, v_i \ra 
\ \geq \ 0 .
\ee
As the matrix in (\ref{eqc2}) is positive semidefinite, it can be
represented as the Grammian of vectors $u_j^r$, which can be chosen to
lie in a Hilbert space of dimension $d_r$ equal to the defect of
$T_r$ (the defect of $T_r$ is the rank of $I -  T_r^\ast T_r$).
So we have
\se\att
\begin{eqnarray}
\label{eqc23}
(1 -\overline{\l_i^1} \l_j^1 ) \la v_j, v_i \ra 
&\=& \la u_j^1, u_i^1 \ra \\
\att\label{eqc24}
(1 -\overline{\l_i^2} \l_j^2 ) \la v_j, v_i \ra 
&\=& \la u_j^2, u_i^2 \ra .
\end{eqnarray}
Multiplying the first equation by $(1 -\overline{\l_i^2} \l_j^2 )$
and the second equation by $(1 -\overline{\l_i^1} \l_j^1 )$, we see
that they are equal. Therefore
\be
\label{eqc3}
(1 -\overline{\l_i^1} \l_j^1 ) \la u_j^2, u_i^2 \ra
\=
(1 -\overline{\l_i^2} \l_j^2 ) \la u_j^1, u_i^1 \ra .
\ee
Reordering equation (\ref{eqc3}), we get
\be
\label{eqc4}
\la u_j^1, u_i^1 \ra \ +\  \overline{\l_i^1} \l_j^1  \la u_j^2, u_i^2
\ra
\=
\la u_j^2, u_i^2 \ra \ +\  \overline{\l_i^2} \l_j^2  \la u_j^1,
u_i^1
\ra.
\ee
Equation~\ref{eqc4} says that there is some unitary matrix 
\be\label{eqc45}
U \= \left( \begin{array}{cc} A & B \\
C & D \end{array} \right) \ : \ \C^{d_1} \oplus \C^{d_2} 
\to \C^{d_1} \oplus \C^{d_2}
\ee
such that
\be
\label{eqc5}
\left( \begin{array}{cc} A & B \\
C & D \end{array} \right)
\left( \begin{array}{c} u_j^1 \\
\l_j^1\ u_j^2 \end{array} \right)
\=
\left( \begin{array}{c}
\l_j^2 u_j^1  \\u_j^2 
\end{array} \right) .
\ee
If the linear span of the vectors $u_j^1 \oplus \l_j^1 u_j^2$ is not
all of $\C^{d_1} \oplus \C^{d_2} $, then $U$ will not be unique. In
this event, we just choose one such $U$.
Define the $d_1 \times d_1$ matrix-valued analytic function $\Psi$ by
\be
\label{eqc55}
\Psi(z) \= A + z B (1 - z D)^{-1} C .
\ee
For any function $\Theta$ of two variables, scalar or
matrix-valued,
define
$$
\Theta^\cup(Z,W) \ := \ \left[ \Theta(Z^\ast, W^\ast) \right]^\ast.
$$
Let $\Phi = \Psi^\cup$, so
$$
\Phi(z) \= A^\ast + z C^\ast (1 - z D^\ast)^{-1} B^\ast .
$$

Equation~\ref{eqc5} implies that
\be
\label{eqc6}
\Psi(\l_j^1)  u_j^1 \=
\left[\Phi (\overline{\l_j^1})\right]^\ast \ u_j^1 \= \l_j^2 u_j^1 .
\ee

Let $s$ be the \sz kernel in the Hardy space $H^2$ of the unit disk, so 
\be
\label{eqc62}
s_\l (z) \= \frac{1}{1 - \overline{\l} z} .
\ee
Let $k_j$ be the vector in $H^2 \otimes \C^{d_1}$ given by
$$
k_j \ := \ s_{\overline{\l_j^1}} \otimes u_j^1 .
$$
Consider the pair of isometries $(M_z, M_\Phi)$
on $H^2 \otimes \C^{d_1}$, where $M_z$ is multiplication by the
coordinate function (times the identity matrix on $\C^{d_1}$) and
$M_\Phi$ is multiplication by the matrix function $\Phi$.
Then
\beq
M_z^\ast  &: & k_j \mapsto \l_j^1 k_j \\
M_\Phi^\ast  &: & k_j \mapsto \l_j^2 k_j .
\eeq
Therefore the map that sends each $v_j$ to $k_j$ gives a unitary
equivalence between $(T_1,T_2)$ and the pair
$(M_z^\ast, M_\Phi^\ast)$ restricted to the span of the vectors $\{
k_j \}_{j=1}^N$. Therefore 
the pair $(M_z^\ast, M_\Phi^\ast)$, acting on the full space $H^2
\otimes \C^{d_1}$, is a co-isometric extension of
$(T_1,T_2)$.

Let $p$ be any polynomial (scalar or matrix valued) in two variables.
We have
\se\att
\begin{eqnarray}
\nonumber
\| p(T_1,T_2 ) \| &\=& \| p(M_z^\ast, M_\Phi^\ast) |_{\vee \{ k_j \} }
 \| \\
\nonumber
&\leq & \| p (M_z^\ast, M_\Phi^\ast ) \|_{H^2 \otimes \C^{d_1}} \\
\nonumber
&= & \| p^\cup (M_z, M_\Phi ) \|_{H^2 \otimes \C^{d_1}} \\
\nonumber
&\leq & \| p^\cup (M_z, M_\Phi ) \|_{L^2 \otimes \C^{d_1} }\\
\label{eqc63}
&= & \| p^\cup  \|_{\partial {V^\cup} } 
\end{eqnarray}
where $V^\cup$ and $V$ are the sets
\se\att
\begin{eqnarray}
\nonumber
V^\cup &\=& \{ (z,w) \in \D^2 \ : \ \det(\Phi(z) - w I) = 0 \} \\
\label{eqc65}
V &\=& \{ (z,w) \in \D^2 \ : \ \det(\Psi(z) - w I) = 0 \} .
\end{eqnarray}
Equality~(\ref{eqc63}) follows from the observation that
\be
\label{eqc64}
\| p^\cup (M_z, M_\Phi ) \|_{L^2 \otimes \C^{d_1} }
\=
\sup_{\theta} \| p^\cup (e^{i\theta} I, \Phi(e^{i\theta}) ) \|,
\ee
where the norm on the right is the operator norm on the $d_1 \times
d_1$ matrices.
Equation~(\ref{eqz2})
shows that, except possibly for the finite set $\sigma(D) \cap \T$,
the matrix $\Phi(e^{i\theta})$ is unitary, and so the norm of any
polynomial applied to $\Phi(e^{i\theta})$ is just the maximum value
of the norm  of the polynomial on 
the spectrum of $\Phi(e^{i\theta})$. By continuity, we obtain 
(\ref{eqc63}).
Taking complex conjugates, (\ref{eqc63}) gives
$$
\| p(T_1,T_2) \| \ \leq \ \| p \|_V,
$$
the desired inequality.

By Theorem~\ref{thmz1},
we see that $V$ and $V^\cup$
are distinguished varieties, and by construction, $V$ contains the
points $\{ (\l_j^1, \l_j^2)  :  1 \leq j \leq N \}$.


\bs
(ii) Now, we drop the assumption that $T = (T_1,T_2)$ be diagonizable.
J. Holbrook proved that
the set of diagonizable commuting matrices is dense in the
set of all commuting matrices \cite{hol82}. So we can assume that there
is a sequence $T^{(n)} = (T_1^{(n)},T_2^{(n)})$ of commuting matrices that converges to 
$T$ in norm and such that each pair satisfies the hypotheses of (i),
\ie each $T^{(n)}$ is a pair of commuting contractions that have
$N$ linearly independent eigenvectors and no unimodular
eigenvalues.
Each $T^{(n)}$ has a unitary $U_n$ associated to it as in (\ref{eqc45}).
By passing to a subsequence if necessary, we can assume that
the defects $d_1$ and $d_2$ are constant, and that 
the matrices $U_n$ converge to a unitary $U$. The corresponding
functions $\Psi_n$ from (\ref{eqc55})
will converge to some function $\Psi$.
Let $q_n(z,w) = \det(\Psi_n(z) -wI)$, and
$q(z,w) = \det(\Psi(z) -wI)$.
Let $V$ be defined by (\ref{eqc65}) for this $\Psi$, and $V_n$ be
the variety corresponding to $\Psi_n$. Notice that the degrees of
$q_n$ are uniformly bounded.

%

Claim: $V$ is non-empty.

Indeed, otherwise it would contain no points of the form $(0,w)$ for
$w \in \D$. That would mean that $\sigma(A) \subseteq \T$, and
so $B$ and $C$ would be zero. That in turn would mean that the
submatrices $A_n$ in $U_n$ would have all their
eigenvalues tending to $\T$, and
hence by (\ref{eqc5}), the eigenvalues of $T_2^{(n)}$ would all tend
to $\T$. Therefore $T_2$ would have a unimodular eigenvalue,
contradicting the hypotheses.

Claim: $V$ is a distinguished variety.
 
This follows from Theorem~\ref{thmz1}.

Claim: Inequality~(\ref{eqc1}) holds.

This follows from continuity.
Indeed, fix some polynomial $p$. For every $\vare >0$, for every $n
\geq n(\vare)$, we have
\beq
\| p(T) \| &\ \leq \ & \vare  +  \| p(T^{(n)}) \| \\
&\leq& \vare  +  \| p \|_{V_n} .
\eeq
We wish to show that 
$$
\lim_{n \to \i} \| p \|_{V_n} \ \leq \ \| p \|_V.
$$
Suppose not. Then there is some sequence $(z_n,w_n)$ in $V_n$
such that 
\be
\label{eqc8}
|p(z_n,w_n) | \ \geq \ \|p \|_V + \vare
\ee
for some $\vare > 0$. Moreover, we can assume that $(z_n,w_n)$
converges to some point $(z_0,w_0)$ in $\overline{\D^2}$.
The point
$(z_0,w_0) $ is in the zero set of $q$, so if it
were in $\D^2$, then it would be in $V$.
Otherwise, $(z_0,w_0)$ must be in $\T^2$.
To ensure that $(z_0,w_0)$ is in $\overline{V}$,
we must rule out the possibility that some sheet of the
zero set of $q$ just grazes the boundary of $\D^2$ without 
ever coming inside.

But this cannot happen. For every $z$ in $\D$, there are 
$d_1$ roots of $\det(\Psi(z) - w I) =0$, and {\it all}
of these occur in $\D$. 
So as $z$ tends to $z_0$ from inside $\D$, one of the $d_1$
branches of 
$w$ must tend to $w_0$ from inside
the disk too.
Therefore $(z_0,w_0)$ is in the closure of $V$, and
(\ref{eqc8}) cannot happen.
\ep

{\bf Remark 1.} 
If $T_1$ has a unimodular eigenvalue $\lambda$, 
then the corresponding eigenspace $\h'$ will be reducing for $T_2$.
Indeed, writing
$$
T_1 \= \left( \begin{array}{cc}
\lambda I & 0 \\
0 & T_1'' \end{array} \right)
\qquad
T_2 \= \left( \begin{array}{cc}
T_2^{\prime} & X \\
0 & T_2'' \end{array} \right),
$$
the commutativity of $T_1$ and $T_2$ means
$X(T_1'' - \lambda) = 0$.
As $\lambda$ is not in the spectrum of $T_1''$, it follows that $X
= 0$.

Therefore
for any polynomial $p$, we have
\be
\label{eqc9}
\| p (T_1, T_2 ) \| \= \max \left( \| p(\lambda I, T_2^\prime) \|,
\| p (T_1'', T_2'' ) \| \right) .
\ee
By von Neumann's inequality for one matrix, the first entry on the
right-hand side of (\ref{eqc9}) is majorized by
$$
\| p \|_{\{ \lambda \times \D \} } .
$$
So if we allow the matrices to have unimodular eigenvalues, we can
still obtain (\ref{eqc1}) by adding to $V$ a finite number of disks
in the boundary of $\D^2$. The new $V$, however, will not be a
distinguished variety.

{\bf Remark 2.} Once one knows And\^o's inequality for matrices, then
it follows for all commuting contractions by approximating them by
matrices --- see \cite{dru83} for an explicit construction.
Of course, the set $V$ must be replaced by the limit points of the sets
that occur at each stage of the approximation, and in general this
may be the whole bidisk.

{\bf Remark 3.} We have actually constructed 
a co-isometric extension of $T$ that is localized to $V$,
and a unitary dilation of $T$ with spectrum contained in $\partial
V$.

\section{The uniqueness variety}
\label{secd}

A {\em solvable Pick problem on $\D^2$} is a set $\{\l_1, \dots,
\l_N \}$ of points in $\D^2$ and a set $\{w_1, \dots, w_N\}$ of
complex numbers such that there is some function $\phi$ of norm less
than or equal to one
in $\hib$ that interpolates (satisfies $\phi(\l_i) = w_i \ \forall
\ 1 \leq i \leq N$).
An {\em extremal Pick problem} is a solvable Pick problem for which
no function of norm less than one interpolates.
The points $ \l_i $ are called the {\em nodes},
and $w_i$ are called the {\em values}. By {\em interpolating
function} we mean any function in the closed
unit ball of $\hib$ that interpolates. 
\bs
Consider the two following examples, in the case $N=2$.

Example 1. Let $\l_1 = (0,0), \l_2 = (1/2,0),
w_1 =0, w_2 = 1/2$. Then a moment's thought reveals that the
interpolating function is unique, and is given by
$\phi(z,w) = z$.

Example 2. Let $\l_1 = (0,0), \l_2 = (1/2,1/2),
w_1 =0, w_2 = 1/2$. Then the interpolating function is far from
unique --- either coordinate function will do, as will any convex
combination of them. (A complete description of all solutions is
given by J.~Ball and T.~Trent in \cite{baltre98}).
But on the distinguished variety $\{(z,z)  :
z  \in  \D \}$, all solutions coincide by Schwarz's lemma.
\bs
For an arbitrary solvable Pick problem, let $\U$ be the set of
points in $\D^2$
on which all the interpolating functions in the closed unit ball of
$\hib$ have the same
value. The preceding examples show that $\U$ may be either the
whole bidisk or a proper subset. In the event that $\U$ is not the
whole bidisk, it is a variety. Indeed, for any $\l_{N+1}$ not in
$U$, there are two distinct values $w_{N+1}$ and $w'_{N+1}$
so that the corresponding $N+1$ point Pick problem has a solution.
By \cite{baltre98,agmc_bid} these problems have interpolating
functions that are rational, of degree bounded by $2(N+1)$. The
set $\U$ must lie in the zero set of the difference of these
rational functions. Taking the intersection over all $\l_{N+1}$ not
in $\U$, one gets that $\U$ is the intersection of the zero sets of
polynomials.
Therefore $\U$ is a variety, and
indeed, by factoring these polynomials into their irreducible
factors, we see that $\U$ is the intersection with the bidisk of
the zero set of one polynomial, together with possibly a finite
number of isolated points. We shall call $\U$ the {\em uniqueness
variety}. (If the problem is not extremal, $\U$ is just the
original set of nodes).

We shall say that an $N$-point extremal Pick problem is {\it minimal} if
none of the $(N-1)$ point subproblems is extremal. The main result of
this section is that if the uniqueness variety is not the whole
bidisk, then it at least contains a distinguished variety running
through the nodes.
If $N=3$, it is shown in \cite{agmc_three} that either $\U = \D^2$
or the minimal extremal problem has a solution that is a function of
one coordinate function only.

\bt
\label{thmd1}
Let $N \geq 2$, and let
$ \l_1, \dots, \l_N $ and $w_1, \dots, w_N$ be the 
data for a minimal extremal Pick problem on the bidisk.
The uniqueness variety $\U$ 
contains a distinguished variety $V$ that contains each of the
nodes.
\et

For a point $\l$ in $\D^2$, we shall write $\l^1$ and $\l^2$ for
the first and second coordinates, respectively. Given a set of
points $\{\l_1, \dots, \l_N\}$ in $\D^2$, an {\em admissible kernel
$K$} is an $N$-by-$N$ positive definite matrix, with all the
diagonal entries $1$, such that
\be
\label{eqd2}
[ ( 1 - \l_i^r \overline{\l_j^r}) K_{ij} ] \ \geq \ 0 \qquad r =
1,2 .
\ee
A theorem of the first author \cite{ag1} 
asserts that a Pick problem on $\D^2$
is solvable if and only if, for every admissible kernel $K$, the
matrix
\be
\label{eqd3}
[(1 - w_i \overline{w_j}) K_{ij}] 
\ee
is positive semi-definite (see \cite{colwer94,baltre98,agmc_bid} 
for alternative proofs).
We shall say that an admissible kernel is {\em active}
if the matrix~(\ref{eqd3}) has a non-trivial null-space, \ie if it
is positive semi-definite but not positive definite.

\begin{lemma}
\label{lemd4}
A solvable Pick problem has an active kernel if and only if it is
extremal.
\end{lemma}
\bp
($\Rightarrow$)
If the problem were not extremal, then for some $\rho < 1$ one
would have
\be
\label{eqd4}
[(\rho^2 I - w_i \overline{w_j}) K_{ij}] \quad \geq \quad 0
\ee
for all admissible kernels. Take $K$ to be an active kernel, 
with $\gamma$ a non-zero vector in the null-space of
$[( I - w_i \overline{w_j}) K_{ij}]$. Then taking the inner product
of the left-hand side of
(\ref{eqd4}) applied to $\gamma$ with $\gamma$ gives $-(1-\rho^2)\|
\gamma \|^2$, which is negative.

($\Leftarrow$) As the problem is extremal, for each $\rho < 1$
there is some admissible kernel $K$ such that $(\rho^2 I - w_i
\overline{w_j}) K_{ij}$ is not positive semi-definite. By
compactness of the set of $N$-by-$N$ positive semi-definite
matrices with $1$'s down the diagonal, there therefore exists some
positive semi-definite $K$, satisfying  (\ref{eqd2}),
and such that (\ref{eqd3}) is not positive definite. 
It just remains to show that this $K$ is actually positive
definite, and therefore a kernel.

Suppose it were not, so for some non-zero vector $v=
(v^1,\dots,v^N)^t$, we have $Kv = 0$.
By (\ref{eqd2}), for each $r=1,2$, the vector $\l^r \cdot v$ (\ie
the vector whose $i^{\rm th}$ component is $\l^r_i v^i$) is also 
in the null-space of $K$. Iterating this observation, one gets that
for any polynomial $p$, the vector $$
p(\l) \cdot v 
\=
\left( \begin{array}{c} p(\l_1) v^1 \\
\vdots\\
p(\l_N ) v^N 
\end{array} \right)
$$
is in the null-space of $K$. Taking $p$ to be a polynomial that is
$1$ at $\l_1$ and zero on the other nodes, we get $K_{11} = 0$, a
contradiction.
\ep

\begin{lemma}
Every admissible kernel on a set $\{ \l_1, \dots, \l_N \}$
can be extended to a continuous admissible kernel 
on a distinguished variety that contains the points $\{ \l_1,
\dots, \l_N \}$.
\label{lemd5}
\end{lemma}
\bp
Let $K$ be an admissible kernel on the set $\{ \l_1, \dots, \l_N
\}$. As it is positive definite, there
are vectors $v_i$  in $\C^N$ such that $K_{ij} = \la v_j, v_i \ra$.
Because $K$ is admissible, Equations~(\ref{eqc23}) and
(\ref{eqc24}) hold. Following the proof of Theorem~\ref{thmc1}, one
gets that for every point $(z,w)$ in the variety $V$ given by
(\ref{eqc65}), one has
non-zero vectors $\widehat{u^1}( z,w)$ and $\widehat{u^2}( z,w)$
such that
$$
\left( \begin{array}{cc} A & B \\
C & D \end{array} \right)
\left( \begin{array}{c} \widehat{u^1}( z,w) \\
z\ \widehat{u^2}( z,w) \end{array} \right)
\=
\left( \begin{array}{c}
w\ \widehat{u^1}( z,w)  \\\widehat{u^2}( z,w)
\end{array} \right) .
$$
Moreover, as the vector $(\widehat{u^1}, \widehat{u^2})^t$
 must just be chosen in the 
null-space of
$$
\left( \begin{array}{cc} A - wI & zB \\
C & zD - I \end{array} \right),
$$
it can be chosen continuously.
When $(z,w)$ is one of the nodes $\l_j$, we choose
\beq
\widehat{u^1} (\l_j^1,\l_j^2) &\=& u_j^1 \\
\widehat{u^2} (\l_j^1,\l_j^2) &\=& u_j^2 .
\eeq

Normalize the vectors so that $$
\| \widehat{u^1} (z,w) \| \= \sqrt{1- |z|^2} .
$$
Now let
$$
k(z,w) \= s_{\bar z} \otimes \widehat{u^1} (z,w) ,
$$
where $s$ is the \sz kernel on the disk as in (\ref{eqc62}).

The desired extension of $K$ to $V$ is given by 
$$
\widehat{K}(\zeta,\lambda) \= \la k(\lambda), k(\zeta) \ra .
$$
This is obviously a kernel that extends  $K$, it is continuous on 
$V \times V$ by construction, and the fact that it
is admissible follows, in the language of Theorem~\ref{thmc1}, from
the fact that $M_z$ and $M_\Psi$  are contractions.
\ep

{\sc Proof of Theorem~\ref{thmd1}:}

(Step 1.) By Lemma~\ref{lemd4}, the problem has an extremal kernel, and
by Lemma~\ref{lemd5}, this kernel can be extended to a distinguished
variety $V$ that contains all the nodes. Let us call the extended
kernel $K$.

Let $\gamma = (\g^1,\dots,\g^N)$ be a non-zero vector in the
null-space
of $[(1 - w_i \overline{w_j}) K_{ij}]$.
Let $\l_{N+1} = (\lN^1,\lN^2)$ be any point in $V$ that is not
one of the original nodes.
Let $\wN$ be some possible value that an interpolating function can
take at $\lN$. As the $(N+1)$ point Pick problem with nodes $\l_1,
\dots,\lN$ and values $w_1,\dots,\wN$ is solvable, and 
as $K$ is admissible, we must have that
$$
[ (1 - w_i \overline{w_j}) K_{ij}]_{i,j=1}^{N+1} \ \geq \ 0.
$$
Therefore, for every $t \in \C$, we have
\be
\label{eqd6}
\la \left[ (1 - w_i \overline{w_j}) K_{ij} \right] 
\left( \begin{array}{c}
\gamma\\
t \end{array} \right) \ , \ 
\left( \begin{array}{c}
\gamma\\
t \end{array} \right)\ra \quad \geq \quad 0.
\ee
As $\gamma$ is in the null-space of $[ (1 - w_i \overline{w_j})
K_{ij}]_{i,j=1}^{N}$, Inequality~(\ref{eqd6}) reduces to 
\be
\label{eqd7}
2 \Re [ \bar t \sum_{j=1}^N (1 - \bar w_j \wN) K_{N+1,j} \g^j]
\ + \ |t|^2 (1 - |\wN|^2) \ \geq \ 0 .
\ee
As this holds for all $t$, we must have that the linear term
vanishes, 
and so we can solve for $\wN$ and get
\se\att
\begin{eqnarray}
\label{eqd799}
\wN \ \left( \sum_{j=1}^N \bar w_j K_{N+1,j} \g^j \right)
&\=&
\sum_{j=1}^N  K_{N+1,j} \g^j .
\\
\att
\wN &\=& \sum_{j=1}^N  K_{N+1,j} \g^j / \sum_{j=1}^N \bar w_j
K_{N+1,j} \g^j .
\label{eqd8}
\end{eqnarray}
As long as both sides of (\ref{eqd799}) do not reduce to zero, 
this gives a formula for $\wN$, which must therefore be unique.
\bs
(Step 2.) So far, we have not used the minimality of the problem.
Minimality ensures that no component of $\gamma$ can be zero, for
otherwise an $(N-1)$ point subproblem would have an active kernel.

Fix one of the nodes, $\l_1$ say, and consider what happens when 
$\lN$ tends to $\l_1$ along some sheet of $V$.
By continuity, $K_{N+1,j}$ tends to $K_{1,j}$ for each $j$.
If $$
\sum_{j=1}^N \bar w_j K_{1,j} \g^j \quad \neq \quad 0,
$$
then by continuity
$$
\sum_{j=1}^N \bar w_j K_{N+1,j} \g^j \quad \neq \quad 0
$$
for $\lN$ in $V$ and close to $\l_1$, and so 
formula~(\ref{eqd8}) gives the unique value that the interpolating
function must take at $\lN$.

Assume instead that 
\be
\label{eqd9}
\sum_{j=1}^N \bar w_j K_{1,j} \g^j \quad = \quad 0 .
\ee
Consider the $N$ point Pick problem with nodes $\l_1, \dots, \l_N$,
and values $w_1+\vare, w_2, \dots , w_N$ for some $\vare$ in
$\C$. If this problem were solvable, then, since $K$ is an
admissible kernel, one would have
\be
\label{eqd10}
[ ( 1 - w_i' \bar w_j') K_{ij}] \quad \geq \quad 0,
\ee
where 
$$
w_i' \= \left\{ \begin{array}{cl}
w_i & i \neq 1 \\
w_1 + \vare & i = 1
\end{array} \right.
$$
Take the inner product of the left-hand side of (\ref{eqd10})
applied to $\gamma$ with $\g$. We get
\se\att
\begin{eqnarray}
\lefteqn{
\sum_{i,j=1}^N ( 1 - w_i' \bar w_j') K_{ij} \g^j \bar \g^i \= } 
\nonumber \\
&&
\sum_{i,j=1}^N ( 1 - w_i \bar w_j) K_{ij} \g^j \bar \g^i
 - 
2 \Re[ \vare \bar \g^1 \sum_{j=1}^N \bar w_j K_{1j} \g^j ]
 -  
|\vare |^2 K_{11} |\g^1|^2
\label{eqd11}
\end{eqnarray}
The first sum in (\ref{eqd11}) vanishes because $\g$ is in the null
space of $[( 1 - w_i \bar w_j) K_{ij}]$. 
The second sum vanishes by hypothesis~(\ref{eqd9}).
Therefore for any $\vare \neq 0$, (\ref{eqd11}) is negative.
This means that the value $w_1$ at $\l_1$ is uniquely determined by
the choice of the other $N-1$ values at $\l_2, \dots, \l_N$. 
Therefore this $(N-1)$ point subproblem must be extremal,
contradicting the minimality hypothesis.

We conclude therefore that (\ref{eqd8}) gives a well-defined
formula for the unique value of $\wN$ at points $\lN$ in $V$ near
the nodes. As we know that some solution to the problem is given by
a rational function, we therefore know that this rational function
gives the unique solution near the nodes. Therefore the union of
the irreducible components of $V$ that contain the nodes is a
distinguished variety contained in $\U$.
\ep

\bs

\quest
Is the distinguished variety constructed in the proof equal to all
of $\U$?
\bs
Given any function on any subset of the bidisk, the result in
\cite{ag1} tells whether it can be extended to a function in the 
closed unit ball of $\hib$. If the set is a distinguished variety,
and the function is analytic on it, is there a better criterion,
which one might think of as solving Problem (b) in the
Introduction?
\quest
How can one tell whether a  function on a distinguished variety 
extends to all of $\D^2$ without increasing its norm?

\end{document}